# Representations of 2-transitive topological groups

Robert A. Bekes

**Abstract**

An analogue of Burnside's Lemma for 2-transitive groups is shown to hold for a class of topological groups. If the group is compact the representation is finite and splits into an irreducible and the constant functions.. If both the group and representation space are noncompact the representation is irreducible.

**Introduction**: Let $G$ be a group acting on a set $X$. The action of $G$ is transitive if for $x$ and $y$ in $X$ there exists $g$ in $G$ such that $gx = y$. The action of $G$ is 2-transitive if for $x_1 \neq x_2$ and $y_1 \neq y_2$ in $X$ there exists $g$ in $G$ such that $gx_1 = y_1$ and $gx_2 = y_2$. Throughout assume that $G$ is a topological group acting on a topological space $X$ such that the map $G \times X \to X$ is continuous and that the action is 2-transitive. If $G$ is finite, using Burnside's Lemma, it can be shown that the associated representation splits into two representations, the identity representation and an irreducible representation orthogonal to the identity representation; see Serre [7], Section 2.3, problem 2.6. For infinite discrete groups Chernoff [2] showed that the associated representation is irreducible. The purpose of this paper is to show that for a class of noncompact $G$ and $X$ the associated representation is irreducible and for compact $G$ and $X$ the representation splits as in the finite case and that $X$ must be finite.

## 1. Noncompact $G$ and $X$:

Let $G$ be a noncompact nondiscrete locally compact and $\sigma$-compact topological transformation group acting faithfully and 2-transitively on a locally compact noncompact not totally disconnected space $X$. Let $x_0$ be in $X$ and $H = \{g \in G \mid gx_0 = x_0\}$, the stabilizer of $x_0$. Then by Kramer [5], Theorems B and 5.14, $G$ is a Lie group and the semidirect product of an abelian normal real vector subgroup $V$, homeomorphic with $\mathbb{R}^n$, with $H$.

The following Lemma is essentially Lipsman [6], pg. 281.

**Lemma 1**: $X$ has a relatively invariant measure.

**Proof:** $G = H \ltimes V$ and so a left Haar measure on $G$ can be written as $dg = dh\,dv$ where $dv$ is a Haar measure on $V$ and $dh$ is a left Haar measure on $H$. Define $\rho: H \to \mathbb{R}$ by $\int_V f(hvh^{-1})dv = \rho(h)\int_V f(v)dv$. Then $\rho$ is a positive continuous function on $H$ with

$\rho(h_1 h_2) = \rho(h_1)\rho(h_2)$ and $\rho(h^{-1}) = \rho(h)^{-1}$.

Since $V$ is abelian $\int_V f(hvv'v^{-1}h^{-1})dv' = \rho(h)\int_V f(vv'v^{-1})dv' = \rho(h)\int_V f(v')dv'$. Therefore $\rho$ can be extended to all of $G$ by $\rho(hv) = \rho(h)$.

Let $f$ be continuous on $G$ with compact support and let $\tilde{f}(gH) = \int_H f(gh)dh$. Then

$$\int_G f(g)\rho(g)^{-1}dg = \int_H\int_V f(hv)\rho(h)^{-1}dvdh = \int_V\int_H f(vh)dhdv$$
$$= \int_{G/H}\int_H f(gh)dhd(gH) = \int_{G/H}\tilde{f}(gH)d(gH).$$

Therefore

$$\int_{G/H}\tilde{f}(g_0 gH)d(gH) = \int_G f(g_0 g)\rho(g)^{-1}dg = \int_G f(g)\rho(g_0^{-1}g)^{-1}dg$$
$$= \rho(g_0)\int_G f(g)\rho(g)^{-1}dg = \rho(g_0)\int_{G/H}\tilde{f}(gH)d(gH).$$

And so $d(gH)$ is relatively invariant on $G/H$. Under the homeomorphism $gH \to gx_0$, $G/H$ is homeomorphic to $X$, and so the result follows. □

**Lemma 2:** Let $x_1 \neq x_2$ in $X$. Then there exists $x_3, x_4, \ldots$ and $\delta > 0$ such that the $x_i$ are all distinct and for $i \neq j$ there exists $g \in G$ with $gx_i = x_1$, $gx_j = x_2$, and $\rho(g) \leq \delta$.

**Proof**: Let $\mu$ be a relatively invariant measure on $X$. Since $G$ is 2-transitive, there exists a $g_0$ such that $g_0 x_2 = x_1$ and $g_0 x_1 = x_2$. Let $\delta = \rho(g_0)$. Suppose we have $x_1, x_2, x_3, \cdots, x_n$ such that for $i \neq j$ there is $g \in G$ with $gx_i = x_1$, $gx_j = x_2$, and $\rho(g) \leq \delta$. Choose $x_{n+1}$ as follows: For each $1 \leq i < n$ let $H_i$ be the stabilizer of $x_i$. Fix $n$. Since $G$ is 2-transitive, $X = \{gx_n \mid g \in H_i\}$. Let $A = \{gx_n \mid \rho(g) \geq 1\}$. Then $A$ has nonempty interior and so $\mu(A) > 0$. Therefore also $\mu(A\setminus\{x_1,\ldots,x_n\}) > 0$. Choose $x_{n+1} \in A\setminus\{x_1,\ldots,x_n\}$. Then for each $1 \leq i < n$ there exists $g_i$ such that $x_{n+1} = g_i x_n$, $g_i x_i = x_i$, and $\rho(g_i) \geq 1$. So $g_i^{-1}x_{n+1} = x_n$, $g_i^{-1}x_i = x_i$, and $\rho(g_i^{-1}) \leq 1$. By the choice of $x_1, x_2, x_3, \cdots, x_n$, there exists $g'$ such that $g'x_n = x_2$, $g'x_i = x_1$, and $\rho(g') \leq \delta$. Then setting $g = g'g_i^{-1}$ we get $gx_{n+1} = x_2$, $gx_i = x_1$, and $\rho(g) = \rho(g'g_i^{-1}) = \rho(g')\rho(g_i^{-1}) \leq \delta$. Therefore the set $x_1, x_2, x_3, \ldots, x_n, x_{n+1}$ has the desired property. □

**Theorem 1**: Let $G$ be as above. The unitary representation associated with the action of $G$ on $X$ is irreducible.

**Proof**: Let $L_2(X,\mu)$ be the Hilbert Space associated with the representation of $G$ on $X$ and let $\langle \, , \, \rangle$ be the inner product. Then $\pi(g)f(x) = \sqrt{\rho(g)^{-1}} f(g^{-1}x)$ defines a unitary representation $\pi$ of $G$ on $L_2(X,\mu)$. See Folland [4], Section 6.1, pg. 154.

Let $U$ be a nonempty set in $X$ with $\mu(U) < \infty$ and let $x_1$ and $x_2$ be such that $x_1 U$ and $x_2 U$ are disjoint. Using Lemma 2 we get distinct $x_1, x_2, x_3, \ldots$ in $X$ and $\delta > 0$ such that for $i \neq j$ there exists $g_{ij} \in G$ with $g_{ij} x_i = x_1$, $g_{ij} x_j = x_2$ and $\rho(g_{ij}) \leq \delta$. Therefore $g_{ij}(x_i U) = (g_{ij} x_i)U = x_1 U$, $g_{ij}(x_j U) = (g_{ij} x_j)U = x_2 U$, and so $\{x_i U\}$ are disjoint.

For any subset $W$ of $X$ let $\xi_W$ denote the characteristic function of $W$. Let $f_n = \sum_{i=1}^n c_i \xi_{x_i U}$ with $c_i \geq 0$. Since $V$ can be identified with $X$ and acts transitively on $X$, there exists $v_i \in V$ such that $v_i x_i = x_1$. By the definition of $\rho$ in Lemma 1, $\rho \equiv 1$ on $V$. Therefore $\mu(x_i U) = \mu(x_1 U)$ and so $\langle f_n, f_n \rangle = \sum_{i=1}^n c_i^2 \mu(x_1 U)$. Let $T$ be a positive intertwining operator for the action $\pi$ of $G$ on $X$. Then $\langle T\xi_{x_i U}, \xi_{x_i U} \rangle = \langle \pi(v_i) T\xi_{x_i U}, \pi(v_i) \xi_{x_i U} \rangle = \langle T\xi_{v_i x_i U}, \xi_{v_i x_i U} \rangle = \langle T\xi_{x_1 U}, \xi_{x_1 U} \rangle$ and for $i \neq j$

$\langle T\xi_{x_i U}, \xi_{x_j U} \rangle = \langle \pi(g_{ij}) T\xi_{x_i U}, \pi(g_{ij}) \xi_{x_j U} \rangle = \rho(g_{ij})^{-1} \langle T\xi_{g_{ij} x_i U}, \xi_{g_{ij} x_j U} \rangle = \rho(g_{ij})^{-1} \langle T\xi_{x_1 U}, \xi_{x_2 U} \rangle$.

Therefore

(1) $\langle Tf_n, f_n \rangle = \sum_{i=1}^n c_i^2 \langle T\xi_{x_i U}, \xi_{x_i U} \rangle + \sum_{i \neq j} c_i c_j \langle T\xi_{x_i U}, \xi_{x_j U} \rangle$

$= \sum_{i=1}^n c_i^2 \langle T\xi_{x_1 U}, \xi_{x_1 U} \rangle + \sum_{i \neq j} c_i c_j \rho^{-1}(g_{ij}) \langle T\xi_{x_1 U}, \xi_{x_2 U} \rangle$.

Since $T$ is positive, $\langle T\xi_{x_1 U}, \xi_{x_2 U} \rangle$ is real.

If $\langle T\xi_{x_1 U}, \xi_{x_2 U} \rangle \geq 0$, from (1) we get

(2) $\langle Tf_n, f_n \rangle \geq \sum_{i \neq j} c_i c_j \rho(g_{ij})^{-1} \langle T\xi_{x_1 U}, \xi_{x_2 U} \rangle$

$\geq \left[ \sum_{i \neq j} c_i c_j \delta^{-1} \right] \langle T\xi_{x_1 U}, \xi_{x_2 U} \rangle$

$= \left[ \left[ \sum_{i=1}^n c_i \right]^2 - \sum_{i=1}^n c_i^2 \right] \delta^{-1} \langle T\xi_{x_1 U}, \xi_{x_2 U} \rangle$.

Now let $c_i = \frac{1}{i}$ and $f = \sum_{i=1}^{\infty} c_i \xi_{x_i U}$. Then since $\sum_{i=1}^{\infty} c_i^2 < \infty$ we have $f \in L_2(X, \mu)$ and $\lim_{n \to \infty} \langle Tf_n, f_n \rangle = \langle Tf, f \rangle < \infty$. Since $\sum_{i=1}^{\infty} c_i = \infty$, letting $n \to \infty$ in (2) we must have $\langle T\xi_{x_1 U}, \xi_{x_2 U} \rangle = 0$.

If $\langle T\xi_{x_1 U}, \xi_{x_2 U} \rangle \leq 0$, from (1) we get

(3) $\langle Tf_n, f_n \rangle = \sum_{i=1}^{n} c_i^2 \langle T\xi_{x_i U}, \xi_{x_i U} \rangle + \sum_{i \neq j} c_i c_j \rho^{-1}(g_{ij}) \langle T\xi_{x_1 U}, \xi_{x_2 U} \rangle$

$\leq \sum_{i=1}^{n} c_i^2 \langle T\xi_{x_i U}, \xi_{x_i U} \rangle + \sum_{i \neq j} c_i c_j \delta^{-1} \langle T\xi_{x_1 U}, \xi_{x_2 U} \rangle$

$= \sum_{i=1}^{n} c_i^2 \langle T\xi_{x_i U}, \xi_{x_i U} \rangle + \left[ \left[ \sum_{i=1}^{n} c_i \right]^2 - \sum_{i=1}^{n} c_i^2 \right] \delta^{-1} \langle T\xi_{x_1 U}, \xi_{x_2 U} \rangle.$

As above, let $c_i = \frac{1}{i}$. Then $\sum_{i=1}^{\infty} c_i^2 < \infty$, $\sum_{i=1}^{\infty} c_i = \infty$, and $\lim_{n \to \infty} \langle Tf_n, f_n \rangle = \langle Tf, f \rangle < \infty$. So with $\langle T\xi_{x_1 U}, \xi_{x_2 U} \rangle \leq 0$, letting $n \to \infty$ in (3) we must have $\langle T\xi_{x_1 U}, \xi_{x_2 U} \rangle = 0$.

Therefore if $\mu(U) < \infty$, $\langle T\xi_{x_1 U}, \xi_{x_2 U} \rangle = 0$ when $x_1 U$ and $x_2 U$ are disjoint.

Since $X$ is homeomorphic with $\mathbb{R}^n$, there is a sequence $\{U_k\}_{k=1}^{\infty}$ of subsets of $X$ with $\mu(U_k) \to 0$ such that each $U_k$ is the disjoint union of $U_{k+1}$ and a translate of $U_{k+1}$ and finite linear combinations of characteristic functions of disjoint translates of $U_k$, $k \geq 1$, are dense in $L_2(X, \mu)$.

Now suppose $W = U \cup xU$ where $\mu(U) < \infty$ and $U$ and $xU$ are disjoint. Then by the above argument, $\langle T\xi_U, \xi_{xU} \rangle = \langle T\xi_{xU}, \xi_U \rangle = 0$. There exists $v$ in $V$ such that $v(xU) = U$. Since $\rho \equiv 1$ on $V$, we also get $\langle T\xi_{xU}, \xi_{xU} \rangle = \langle \pi(v)T\xi_{xU}, \pi(v)\xi_{xU} \rangle = \langle T\xi_{v(xU)}, \xi_{v(xU)} \rangle = \langle T\xi_U, \xi_U \rangle$. Let $\lambda = \frac{\langle T\xi_W, \xi_W \rangle}{\mu(W)}$. Then $\lambda = \frac{\langle T\xi_U, \xi_U \rangle + \langle T\xi_{xU}, \xi_{xU} \rangle}{2\mu(U)} = \frac{\langle T\xi_U, \xi_U \rangle}{\mu(U)}$. So for any such decomposition, $\lambda$ is independent of $U$ and $W$ and so $\langle T\xi_W, \xi_W \rangle = \lambda \langle \xi_W, \xi_W \rangle$ and $\langle T\xi_U, \xi_U \rangle = \lambda \langle \xi_U, \xi_U \rangle$. Therefore $\langle T\xi_{U_k}, \xi_{U_k} \rangle = \lambda \langle \xi_{U_k}, \xi_{U_k} \rangle$ for all $k$ and so $T = \lambda I$. It then follows that the representation $\pi$ is irreducible. □

**The *ax+b* group:**

Let $G$ be the $ax+b$ group acting on $\mathbb{R}$ by $x \mapsto ax+b$ where $a \neq 0$. If $x_1 \neq x_2$ and $y_1 \neq y_2$, the system $\begin{bmatrix} x_1 & 1 \\ x_2 & 1 \end{bmatrix} \begin{bmatrix} a \\ b \end{bmatrix} = \begin{bmatrix} y_1 \\ y_2 \end{bmatrix}$ has a solution. Therefore the action of $G$ is 2-transitive on $\mathbb{R}$ and so the representation of $G$ on $\mathbb{R}$ is irreducible. See Conrad [3], example 4.3.

This result also follows from the representation theory of semidirect products, see Folland [4], Section 6.7, pg. 189.

**Discrete groups:**

The argument in Theorem 1 simplifies for infinite discrete groups acting on an infinite discrete set $X$. To prove irreducibility, start by selecting distinct $x_1, x_2, \ldots$ in $X$ and replacing $x_i U$ by the singleton $\{x_i\}$ and $W$ with $\{x_k, x_1\}$. This case is proved in Chernoff [2].

**Example**:

Let $G$ be the group of permutations on the integers that move only a finite number of integers. Then $G$ acts 2-transitively and so the associated representation is irreducible.

**2. Compact $G$ and $X$:**

In this section $G$ is a compact topological group acting 2-transitively on a nontrivial compact topological space $X$. Let $\mu$ be a $G$ invariant measure on $X$. Let $L_2(X, \mu)$ be the Hilbert space of the associated representation of $G$ on $X$ and let $\langle\ ,\ \rangle$ be the inner product. The representation of $G$ on $L_2(X, \mu)$ is defined by $\pi(g)f(x) = f(g^{-1}x)$. Set $\mu(X) = 1$.

The idea for the following proof was communicated to me by Don Blasius. It is an extension of the argument in Serre [7], Section 2.3, problem 2.6 to compact groups.

**Theorem 2**: $X$ is finite and the representation $\pi$ splits into an irreducible representation and projection onto the space of constant functions.

**Proof** : Let $\chi$ be the character of the representation $\pi$ of $G$ on $X$. Let $\bar{\pi}$ be the conjugate representation of $\pi$. Then $\bar{\chi}$ is the character of $\bar{\pi}$ and so the character of $\pi \otimes \bar{\pi}$ is $|\chi|^2$. Let 1 be the character of the unit representation of $G$. In $L_2(X, \mu)$, $\langle |\chi|^2, 1 \rangle$ equals the number of times the unit representation is contained in the representation $\pi \otimes \bar{\pi}$ of $G$ on $X \times X$. Because $G$ is 2-transitive there are 2 orbits of $G$ in $X \times X$, the diagonal and the off diagonal. Because of transitivity the measure $\mu$ is either purely discrete or continuous. Suppose that $\mu$ is continuous. Then the diagonal in $(X \times X, \mu \times \mu)$ has measure zero. Since $\langle |\chi|^2, 1 \rangle$ can be interpreted as is the number of orthogonal one dimensional $G$ invariant subspaces in $L_2(X \times X, \mu \times \mu)$. Such subspaces are determined by a single function supported on a $G$ orbit in $X \times X$. Therefore $\langle |\chi|^2, 1 \rangle = 1$ when $\mu$ is continuous. But $\langle |\chi|^2, 1 \rangle = \langle \chi, \chi \rangle =$ the number of irreducible subrepresentations of $G$ in $\pi$. The space of constant functions on $G$ is invariant under $G$ and so is its orthogonal complement. Since $G$ is compact, $\pi$ restricted to the orthogonal complement of the space of constant functions on $G$ has an irreducible subrepresentation. So if $\mu$ is continuous we must have $\langle |\chi|^2, 1 \rangle = \langle \chi, \chi \rangle > 1$. Therefore $\mu$ can't be continuous and hence must be purely discrete. Therefore $X$ must be finite. It follows that $\langle |\chi|^2, 1 \rangle = \langle \chi, \chi \rangle = 2$ and so $\pi$ splits into an irreducible representation and projection onto the constant functions. □

## 3. Noncompact $G$ and compact $X$:

In this case if $G$ acts 2-transitively on $X$ and there is a $G$ invariant measure on $X$ the argument in Theorem 2 holds and the representation splits as the sum of an irreducible and the projection onto the constant functions.

If there is no $G$ invariant measure on $X$ the situation is more complicated. Let $G = SL(2, \mathbb{R})$ and let $X = \mathbb{RP}^1$, the real projective line. Then it is shown in Conrad [3], Theorem 4.21, that the action of $G$ on $X$ is 2-transitive. It follows from Casselman [1] page 16, that $X \cong G/B$ where $B$ is the Borel subgroup of $G$ and the representation on $X$ is, via normalized induction, $\text{Ind}_B^G \delta_B^{-1/2}$ where $\delta \begin{bmatrix} t & x \\ 0 & t^{-1} \end{bmatrix} = t^2$. By Casselman [1] Proposition 8.7 with $s = -1$, $m = 1$, and $n = 0$, the orthogonal complement of the projection onto the space of constant functions on $X$ splits into two infinite dimensional subrepresentations.

**Acknowledgements:** I would like to thank Kenneth A. Ross for his suggestions and critique of this paper and Stephen DeBacker for the reference to Casselman [1].

# REFERENCES

[1] Bill Casselman, Analysis on SL(2) representations of $SL_2(\mathbb{R})$, https://pdfs.semanticscholar.org/3556/1bd7456d3fa2892d087e497472109a993f1f.pdf?_ga=2.268481145.921800545.1528385167-1259800467.1528385167

[2] Paul R. Chernoff, Irreducible representations of infinite-dimensional transformation groups and Lie algebras I, Journal of Functional Analysis, 130, 255-282 (1995).

[3] Keith Conrad, Transitive group actions, http://www.math.uconn.edu/~kconrad/blurbs/grouptheory/transitive.pdf

[4] Gerald B. Folland, *A course in harmonic analysis*, CRC Press, Boca Raton, Florida, 1995.

[5] Linus Kramer, Two-transitive Lie groups, Journal für die reine und angewandte Mathematik, 563 (2003), 83—113.

[6] Ronald L. Lipsman, The Penny-Fujiwara Plancherel formula for abelian symmetric spaces and completely solvable homogeneous spaces, Pacific Journal of Mathematics, Vol 151, No 2, 1991.

[7] J.-P. Serre, *Linear representations of finite groups*, Springer-Verlag, New York, 1977.

Robert A. Bekes
Mathematics and Computer Science Dept.
Santa Clara University
Santa Clara, CA 95053
rbekes@scu.edu